\documentclass[1p]{elsarticle}
\usepackage{amssymb}
\usepackage{amsfonts}

\newtheorem{theorem}{Theorem}
\newtheorem{conjecture}[theorem]{Conjecture}

\newproof{pf}{Proof}

\begin{document}
\title{A note on asymptotically optimal neighbour sum distinguishing colourings}

\author{Jakub Przyby{\l}o\fnref{grantJP,MNiSW}}
\ead{jakubprz@agh.edu.pl, phone: 048-12-617-46-38,  fax: 048-12-617-31-65}

\fntext[grantJP]{Supported by the National Science Centre, Poland, grant no. 2014/13/B/ST1/01855.}
\fntext[MNiSW]{Partly supported by the Polish Ministry of Science and Higher Education.}

\address{AGH University of Science and Technology, al. A. Mickiewicza 30, 30-059 Krakow, Poland}

\begin{abstract}
The least $k$ admitting a proper edge colouring $c:E\to\{1,2,\ldots,k\}$ of a graph $G=(V,E)$ without isolated edges
such that $\sum_{e\ni u}c(e)\neq \sum_{e\ni v}c(e)$ for every $uv\in E$ is denoted by $\chi'_{\Sigma}(G)$.
It has been conjectured that $\chi'_{\Sigma}(G)\leq \Delta + 2$ for every connected graph of order at least three different from the cycle $C_5$, where $\Delta$ is the maximum degree of $G$. It is known that $\chi'_{\Sigma}(G) = \Delta + O(\Delta^\frac{5}{6}\ln^\frac{1}{6}\Delta)$ for a graph $G$ without isolated edges. We improve this upper bound to $\chi'_{\Sigma}(G) = \Delta + O(\Delta^\frac{1}{2})$ using a simpler approach involving a combinatorial algorithm enhanced by the probabilistic method.
The same upper bound is provided for the total version of this problem as well.
\end{abstract}


\maketitle

\section{Introduction}
Every simple graph of order at least two
contains a pair of vertices of the same degree.
In this sense, there are no \emph{irregular graphs}, apart from the 1-vertex exception.
A possible alternative definitions of these were studied e.g.
by Chartrand, Erd\H{o}s and Oellermann in~\cite{ChartrandErdosOellermann}, with no definite solution.
Thus, Chartrand et al.~\cite{Chartrand} proposed
a graph invariant designed to measure the level of irregularity of a graph instead.
Suppose that given a graph $G=(V,E)$ we wish to construct a multigraph
with pairwise distinct vertex degrees of it by multiplying some edges of $G$.
The least $k$ so that we are able to do this
using at most $k$ copies of every edge (counting the original one in)
is known as the \emph{irregularity strength} of $G$ and denoted by $s(G)$, see~\cite{Chartrand}.
This can be defined in terms of labellings or colourings as well.
Consider a (not necessarily proper) edge colouring
$c:E\to\{1,2,\ldots,k\}$ and denote the \emph{weighted degree} of a vertex $v\in V$ as
$d_c(v):=\sum_{u\in N(v)}c(uv).$
This shall also be simply called the \emph{sum at} or \emph{of} $v$.
Then $s(G)$ is equal to the least $k$ admitting such a colouring $c$ with $d_c(u)\neq d_c(v)$ for every $u,v\in V$ with $u\neq v$.
It is clear that $s(G)$ is well defined for all graphs
containing no isolated edges and at most one isolated vertex.
The irregularity strength was studied in numerous papers, e.g.~\cite{Aigner,Bohman_Kravitz,Lazebnik,Faudree,Frieze,KalKarPf,Lehel,MajerskiPrzybylo2,Nierhoff,Przybylo,irreg_str2},
and gave rise to now fast-developing discipline, which might be referred to as \emph{additive graph labelings},
or more generally -- \emph{vertex distinguishing graph colourings}.

An intriguing successor of the problem above is the so-called \emph{1--2--3 Conjecture}, where the only difference is that we require so that $d_c(u)\neq d_c(v)$ only for the \emph{neighbours}, i.e. when $uv$ forms an edge in $G$. For such a goal,
Karo\'nski \L uczak and Thomason~\cite{123KLT} conjectured
 the colours $1,2,3$
 to be sufficient for every graph without isolated edges, see also e.g.~\cite{Louigi30,Louigi}.
So far the colours $1,2,3,4,5$ are known to suffice, see~\cite{KalKarPf_123}.
Another important problem of this field concerns
proper edge colourings where every vertex is
required to differ from its neighbours with respect to the set of its incident colours.
The least number of colours needed to construct such a (proper) edge colouring of a graph $G$ is called its \emph{neighbour set distinguishing index} or \emph{adjacent strong chromatic index}, and denoted by $\chi'_a(G)$, see~\cite{Zhang} by Zhang, Liu and Wang, where it was conjectured that $\chi'_a(G)\leq \Delta(G)+2$ for every connected graph of order at least three, unless $G$ is the cycle $C_5$. (Note that due to the properness of the colourings investigated, we obviously have $\chi'_a(G)\geq \Delta(G)$ or $\chi'_a(G)\geq \Delta(G)+1$ depending on the class of $G$, cf. Vizing's Theorem.)
This conjecture, commonly referred to as the Zhang's Conjecture, has been confirmed for several graph classes, see e.g.~\cite{Akbari,BalGLS,HocqMont,Hornak_planar,WangWang}.
In general Hatami~\cite{Hatami} proved that for any graph $G$ with no isolated edges and with maximum degree $\Delta>10^{20}$, $\chi'_a(G)\leq \Delta+300$, using a multistage probabilistic argument.

We shall focus on a hybrid of the two problems above, where we investigate proper edge colourings, but require of these a much stronger property than distinguishing neighbours by sets, namely
the same one as considered within the 1--2--3 Conjecture. A proper colouring $c:E\to\{1,2,\ldots,k\}$ of the edges of a graph $G=(V,E)$ is called \emph{neighbour sum distinguishing} if $d_c(u)\neq d_c(v)$ for every $uv\in E$. The least $k$ admitting such a colouring is denoted by $\chi'_\Sigma(G)$ and called the \emph{neighbour sum distinguishing index} of $G$. Though such a condition is indeed significantly
stronger than the one investigated within the Zhang's Conjecture, what in particular implies that $\chi'_a(G)\leq \chi'_\Sigma(G)$ for every graph $G$ without isolated edges, the following strengthening of this conjecture was posed in~\cite{FlandrinMPSW}.
\begin{conjecture}[{\cite{FlandrinMPSW}}]\label{Flandrin_et_al_Conjecture}
If $G$ is a connected graph of order
at least three different from the cycle $C_5$, then
$\chi'_{\Sigma}(G) \leq \Delta(G) + 2$.
\end{conjecture}
In fact we are not aware of any graph for which the values of the two parameters differ, see \cite{BonamyPrzybylo,DongWang_mad,FlandrinMPSW,HuChenLuoMiao,WangChenWang_planar,YuQuWangWang} for partial results concerning special graph classes. In general it is known that $\chi'_{\Sigma}(G)\leq 2\Delta(G)+{\rm col}(G)-1\leq 3\Delta(G)$ \cite{Przybylo_CN_1} and $\chi'_{\Sigma}(G)\leq \Delta(G)+3{\rm col}(G)-4$ \cite{Przybylo_CN_2}, where ${\rm col}(G)$ is the colouring number of a graph $G$ (containing no isolated edges).
Unfortunately, the approach of Hatami from~\cite{Hatami} could not be adapted to this more demanding setting. Instead, the following asymptotic confirmation of Conjecture~\ref{Flandrin_et_al_Conjecture}, implying that $\chi'_\Sigma(G)\leq (1+o(1))\Delta(G)$ was provided in~\cite{Przybylo_asym_optim} by means of a different probabilistic approach.
\begin{theorem}\label{main_JP_Th_as}
There is a constant $\Delta_0$ such that
if $G$ is a graph without isolated edges and with maximum degree $\Delta\geq \Delta_0$, then
\begin{equation}\label{main_JP_ineq}
\chi'_{\Sigma}(G) \leq \Delta+50\Delta^{\frac{5}{6}}\ln^\frac{1}{6}\Delta.
\end{equation}
\end{theorem}
The aim of this paper is to present a new simpler approach towards proving a stronger upper bound of the form $\chi'_{\Sigma}(G) \leq \Delta+O(\sqrt{\Delta})$. Our technique yields also the same upper bound for the total correspondent of $\chi'_{\Sigma}$, cf. Theorem~\ref{th_new_asymptotic_edge_and_total} below,
where we investigate the least $k$ so that there exists a proper total colouring $c:V\cup E\to\{1,2,\ldots,k\}$ of a graph $G=(V,E)$
such that $c(u)+d_c(u)\neq c(v)+d_c(v)$ for every edge $uv\in E$, and denote it by $\chi''_{\Sigma}(G)$.
For a given vertex $v$, the sum $c(v)+d_c(v)$ shall be called its \emph{total sum} and denoted by $d^t_c(v)$.
In~\cite{PilsniakWozniak_total}
it was conjectured that $\chi''_{\Sigma}(G)\leq \Delta(G)+3$ for every graph $G$, see~\cite{total_sum_planar,LiLiuWang_sum_total_K_4,PilsniakWozniak_total} for partial results concerning this.
Note that proving this apparently very challenging conjecture would require improving the following best known upper bound due to Molloy and Reed
concerning the well known Total Colouring Conjecture (that $\chi''(G)\leq \Delta(G)+2$ for every $G$) posed by Vizing~\cite{Vizing2} and independently by Behzad~\cite{Behzad}.
\begin{theorem}[\cite{MolloyReedTotal}]\label{MolloyReedTh}
There exists a constant $\Delta_1$ such that for every graph $G$ with $\Delta(G)\geq\Delta_1$, $\chi''(G)\leq\Delta(G)+10^{26}$.
\end{theorem}
In general it has been known thus far that $\chi''_{\Sigma}(G)\leq \Delta(G)+ \lceil\frac{5}{3}{\rm col}(G)\rceil$, \cite{Przybylo_CN_3}, and analogously as in the case of the edge protoplast, $\chi''_{\Sigma}(G) \leq \Delta+O(\Delta^{\frac{5}{6}}\ln^\frac{1}{6}\Delta)$ for any graph $G$, see~\cite{Przybylo_asym_optim_total}.
It is worth mentioning that also a total correspondent of the Zhang's Conjecture is known, see e.g.~\cite{CokerJohanson,Zhang_total}.

\section{Main Result and Tools}

Our main result is the following.
\begin{theorem}\label{th_new_asymptotic_edge_and_total}
There is a constant $\Delta_2$ such that for every graph $G$ with maximum degree $\Delta\geq \Delta_2$,
$\chi''_\Sigma(G)\leq \Delta+95\sqrt{\Delta}$ and $\chi'_\Sigma(G)\leq \Delta+95\sqrt{\Delta}$ if $G$ has no isolated edges.
\end{theorem}
Note that this implies that $\chi''_\Sigma(G)= \Delta+O(\sqrt{\Delta})$ for
all graphs and $\chi'_\Sigma(G)= \Delta+O(\sqrt{\Delta})$ for
graphs without isolated edges.
Within the proof we shall several times make use of two classical
tools of the probabilistic method,
the Lov\'asz Local Lemma, see e.g.~\cite{AlonSpencer}, and the Chernoff Bound, see e.g.~\cite{JansonLuczakRucinski}.
\begin{theorem}[\textbf{The Local Lemma; Symmetric case}]\label{LLL-symmetric}
Let $A_1,A_2,\ldots,A_n$ be events in an arbitrary probability space.
Suppose that each event $A_i$ is mutually independent of a set of all the other
events $A_j$ but at most $D$, and that ${\rm \emph{\textbf{Pr}}}(A_i)\leq p$ for all $1\leq i \leq n$. If
$$ ep(D+1) \leq 1,$$
then $ {\rm \emph{\textbf{Pr}}}\left(\bigcap_{i=1}^n\overline{A_i}\right)>0$.
\end{theorem}
\begin{theorem}[\textbf{Chernoff Bound}]\label{ChernofBoundTh}
For any $0\leq t\leq np$:
$${\rm\emph{\textbf{Pr}}}({\rm BIN}(n,p)>np+t)<e^{-\frac{t^2}{3np}}~~{and}~~{\rm\emph{\textbf{Pr}}}({\rm BIN}(n,p)<np-t)<e^{-\frac{t^2}{2np}}\leq e^{-\frac{t^2}{3np}}$$
where ${\rm BIN}(n,p)$ is the sum of $n$ independent Bernoulli variables, each equal to $1$ with probability $p$ and $0$ otherwise.
\end{theorem}
Note that if $X$ is a sum of $n\leq k$ (where $k$ does not have to be an integer)
random independent Bernoulli variables, each equal to $1$ with probability at most $q\leq 1$,
then $\mathbf{Pr}(X>kq+t)< e^{-\frac{t^2}{3kq}}$ for $0\leq t\leq\lfloor k\rfloor q$.
(It is sufficient to consider a variable $Y$ with binomial distribution ${\rm BIN}(\lfloor k\rfloor,q)$.)

\section{Proof of Theorem~\ref{th_new_asymptotic_edge_and_total}}

It is sufficient to consider a graph $G=(V,E)$ without isolated edges (and vertices) and with sufficiently large maximum degree $\Delta$, i.e. greater than some explicit constant $\Delta_2$ -- large enough so that all inequalities below hold.
We shall thus not specify $\Delta_2$, though we
assume that $\Delta_2\geq \Delta_1$ from Theorem~\ref{MolloyReedTh} above.

We shall consider edge and total colourings simultaneously.
We begin by fixing an arbitrary proper total colouring $c:V\cup E\to\{1,\ldots,\Delta+K\}$ with $K=10^{26}$
following from Theorem~\ref{MolloyReedTh}.
In case of edge colourings exclusively we then set $c(v)=0$ for every $v\in V$ and still call such a total colouring \emph{proper}.

\subsection{Preprocessing Small Degree Vertices}

Let
$$V_S=\{v\in V:d_G(v)\leq \Delta/4\}$$
be the set of vertices with relatively small degrees,
and set $G_S=G[V_S]$. Let $E_I$ be the set of isolated edges in $G_S$. Note that every vertex in $V$
can be incident with at most one edge from $E_I$ in $G$, and analogously every edge in $E$ may be adjacent with at most one element of $E_I$.
Let
$$V_M=\{v\in V: \sqrt{\Delta}+1\leq d_G(v)\leq \Delta/4\}.$$
For every $v\in V_M$ we randomly, independently and uniformly
choose one of its incident edges in $G-E_I$ (hence the probability that any such edge is chosen by $v$ equals at most $1/\sqrt{\Delta}$).
Denote the auxiliary graph induced by the chosen edges by $G_A$.
For every
$v\in V$, by estimating the number of edges chosen to $G_A$ by the neighbours of $v$ (and keeping in mind the one possibly chosen by $v$),
by the Chernoff Bound (and due to the fact that $d_G(v)\leq \Delta$),
$$\mathbf{Pr}(d_{G_A}(v)-1 > \sqrt{\Delta}+\sqrt[3]{\Delta}) < e^{-(\sqrt[3]{\Delta})^2/3\sqrt{\Delta}}<\Delta^{-3}.$$
If we now denote by $A_v$ the event that $d_{G_A}(v) > \sqrt{\Delta}+\sqrt[3]{\Delta}+1$, then it is mutually independent of all other events of this type except those $A_u$ with $u$ at distance at most $2$ from $v$, i.e. at most $\Delta^2+1$ of these.
Hence by the Lov\'asz Local Lemma, with positive probability, $d_{G_A}(v)\leq (1+o(1))\sqrt{\Delta}$ for every $v\in V$. We fix any $G_A$ with this feature.

We then greedily, one by one, modify the colours of the edges of $G_A$ using new at most $2\sqrt{\Delta}(1+o(1))$ colours being the least yet unused positive integers, so that the total colouring is still proper and there are no sum conflicts between the ends of any edge in $E_I$ with at least one end of degree at least $\sqrt{\Delta}+1$. Then, subsequently for every remaining edge in $E_I$, we choose one of its adjacent edges and modify its colour if necessary, using some colour in the range $\{1,\ldots,\Delta-1+\lceil\sqrt{\Delta}\rceil+2+1+1\}$ to get read of potential conflicts between the ends of these edges from $E_I$ (keeping the total colouring proper).
The obtained proper total colouring we denote by $c'$.
Set $m'=\max_{z\in V\cup E}c'(z)$ and note that
at this point,
\begin{equation}\label{m'upperbound}
m'\leq\Delta+2\sqrt{\Delta}(1+o(1)).
\end{equation}

\subsection{Base Subgraph $H$}

Now independently for every edge incident with some vertex in
$$V_L=\{v\in V: d_G(v)\geq \Delta/32\}$$
in $G$
we pick it with probability $6/\sqrt{\Delta}$
and denote the graph induced by all the picked edges by $H$.
For every vertex $v\in V$ we denote by $A'_v$ the event that $v\in V\smallsetminus V_L$ and $d_H(v) > 6\sqrt{\Delta}/32 +\sqrt[3]{\Delta}$,
while by $A''_v$ we denote the event that $v\in V_L$ and $|d_H(v)-6d_G(v)/\sqrt{\Delta}|>\sqrt[3]{\Delta}$.
Similarly as above we shall show that with positive probability none of these events occurs. For any $v\in V$, by the Chernoff Bound
(as $d_G(v) < \Delta/32$ if $v\in V\smallsetminus V_L$),
$$\mathbf{Pr}(A'_v) < e^{-(\sqrt[3]{\Delta})^2/(9\sqrt{\Delta}/16)}<\Delta^{-2}$$
and analogously
$$\mathbf{Pr}(A''_v) < 2e^{-(\sqrt[3]{\Delta})^2/(18d_G(v)/\sqrt{\Delta})}<\Delta^{-2}$$
(as $d_G(v)\leq \Delta$).
Since every event $A'_v$ and $A''_v$ is mutually independent of all other events of these types except those corresponding to vertices at distance at most $1$ from $v$, i.e., at most $2\Delta+1$ other events, the Local Lemma implies that with positive probability none of the events $A'_v,A''_v$ with $v\in V$ occurs. We fix any $H$ with this feature and call it the \emph{base subgraph} of $G$.
Then $d_H(v) \leq 3\sqrt{\Delta}/16 +\sqrt[3]{\Delta}$ for every $v\in V\smallsetminus V_L$ and $|d_H(v)-6d_G(v)/\sqrt{\Delta}| \leq \sqrt[3]{\Delta}$ for $v\in V_L$, hence:
\begin{enumerate}
\item[(I)] $d_H(v)\leq 6\sqrt{\Delta}+\sqrt[3]{\Delta}$ for $v\in V$;
\item[(II)] $d_H(v)\geq 3\sqrt{\Delta}/16 - \sqrt[3]{\Delta}$ for $v\in V_L$.
\end{enumerate}
We shall further introduce colour modifications almost exclusively on the edges of the base subgraph $H$ (and of $G_S$ at the very end).

\subsection{Algorithm Distinguishing Large Degree Vertices}

We shall now use an algorithm partly inspired by~\cite{KalKarPf} in order to distinguish sums of al the neighbours in $V_L$.
At every step of the algorithm below, by $c_T(z)$ we shall understand the contemporary colour of any given edge or vertex $z$ of $G$
(hence $d^t_{c_T}(v)$ shall denote the contemporary sum of a vertex $v$).
Let
$$B=\left\lceil 23\sqrt{\Delta}+3\sqrt[3]{\Delta} \right\rceil.$$

First one by one we greedily recolour the edges in $H$ properly modulo $B$ (i.e. so that the new colours of any two edges adjacent in $H$ are not congruent modulo $B$) with integers in $[m'+B+1,m'+2B]$
so that there are no sum conflicts modulo $B$
between the ends of any edge in $E_I$ adjacent with at least one edge of $H$ (what is feasible due to (I)).
As we have used new colours for this purpose, the obtained total colouring of $G$ is still proper.

Let $H_1,\ldots,H_p$ be the components of $H$. We analyze these
one after another in the stated order.
In each of these we subsequently order the vertices into a sequence so that every vertex except the last one has a neighbour with bigger index in this sequence -- so called \emph{forward neighbour} in $H$; analogously we define a \emph{backward neighbour} of a vertex in $H$. The concatenation of such consecutive sequences for $H_1,\ldots,H_p$ yields an ordering $\pi$ of all the vertices of $H$.

Consider $H_q$, $1\leq q\leq p$, and suppose $H_1,\ldots,H_{q-1}$ has already been analyzed.
Let $v_1,\ldots,v_n$ be the consecutive (all) vertices of $H_q$ in the sequence $\pi$.
We analyze these one after another,
and the moment a given vertex $v\in V(H_q)\cap V_L$ is analyzed, we admit (if necessary):
\begin{enumerate}
\item[(i)] adding or subtracting $B$ to colours of its \emph{backward edges} in $H$ (i.e. joining $v$ with its backward neighbours in $H$);
\item[(ii)] adding $B$ to colours of its \emph{forward edges} in $H$ (i.e. joining $v$ with its forward neighbours in $H$), with one exception, when $e$ is the last forward edge of $v$ (joining $v$ with its forward neighbour $v_l$ in $H$ with maximal $l$), when we admit adding any integer from the set $\{0,1,\ldots,B-1\}$ to the colour of $e$ so that the colouring of $H$ remains proper modulo $B$.
\end{enumerate}
(We introduce no changes if $v\notin V_L$.)
Note that the colour of every edge may be modified at most twice within our algorithm, hence at its end the colouring shall be proper, as we shall use new colours in the range $[m'+1,m'+4B]$ for the edges of $H$ and by~(\ref{m'upperbound}) the largest colour used shall be
bounded from above by
\begin{equation}\label{maximal_colour_ineq}
\Delta+(2+o(1))\sqrt{\Delta}+4B\leq \Delta+94\sqrt{\Delta}+o(\sqrt{\Delta}) \leq \Delta+95\sqrt{\Delta}.
\end{equation}

The moment a given vertex $v\in V(H_q)\cap V_L$ is being analyzed, we associate it with a two-element set $P_v$ of integers from the family
$$\mathcal{P}:=\{\{2iB+j,(2i+1)B+j\}:i\in\mathbb{Z}, j\in\{0,1,\ldots,B-1\}\}$$
making sure that ever since such assignment, $d^t_{c_T}(v)\in P_v$ and that $P_v\cap P_u=\emptyset$ for every neighbour $u\in V_L$ of $v$ in $G$ that has its set $P_u$ already prescribed (hence preceding $v$ in the sequence $\pi$). Only vertices $v_{n-1}$ and $v_n$ shall be handled with in a special manner. This shall assure that at the end of this algorithm $d_{c_T}(u)\neq d_{c_T}(v)$ for every $uv\in E$ with $u,v\in V_L$.

We shall additionally be careful while introducing the colour changes so that:
\begin{enumerate}
\item[(a)] the edge colouring of $H$ remained proper modulo $B$ (what was already mentioned above);
\item[(b)] if $e$
is the last forward edge of some vertex and $e$ is adjacent with an edge $f\in E_I$, then while modifying the colour of $e$ we make sure that the ends of $f$ are sum-distinguished modulo $B$.
\end{enumerate}

Suppose we are about to analyze $v\in V(H_q)\cap V_L\smallsetminus\{v_{n-1},v_n\}$ and thus far all the conditions and requirement listed above have been fulfilled.
Let $e$ be the last forward edge of $v$ in $H$.
Note that as $d_H(v)\leq 6\sqrt{\Delta}+\sqrt[3]{\Delta}$ by (I),
then due to (ii) we have at least $B$ options for the colour of $e$, at most $2(6\sqrt{\Delta}+\sqrt[3]{\Delta}-1)+1$ of which might be forbidden due to requirements (a) and (b). The colour of every remaining edge incident with $v$ in $H$ can on the other hand be modified (if necessary) by $B$ consistently with our requirements (i.e. $B$ can either be added or subtracted from the colour of each such edge $uv$ by (i) or (ii) so that $u\in P_u$ if $u$ has $P_u$ assigned). By (II) we thus have
at least $(3\sqrt{\Delta}/16 - \sqrt[3]{\Delta})(B-12\sqrt{\Delta}-2\sqrt[3]{\Delta}+1)>2\Delta$ distinct sums available for $v$ attainable via admissible changes on the edges incident with $v$. These contain elements of more than $\Delta$ pairs from $\mathcal{P}$.
We choose one of these pairs which is disjoint with the pairs already associated with neighbours of $v$ in $G$, denote it by $P_v$ and perform
admissible changes on the edges incident with $v$ so that $d^t_{c_T}(v)\in P_v$ afterwards.

Assume now that we have analyzed as above all vertices in $V(H_q)\cap V_L\smallsetminus\{v_{n-1},v_n\}$. We may assume that $\{v_{n-1},v_n\}\cap V_L\neq \emptyset$ (as there is nothing to do otherwise).
Suppose first that $\{v_{n-1},v_n\}\subseteq V_L$.
Then we first add to the colour of $v_{n-1}v_n$ some integer in $\{0,1,\ldots,B-1\}$ consistent with (a) and (b), hence analogously as above we have at least $B-12\sqrt{\Delta}-2\sqrt[3]{\Delta}+1$ available options (distinct modulo $B$) for the resulting colour of $v_{n-1}v_n$.
We choose it so that afterwards each of $v_{n-1}$ and $v_n$ has at most $2\sqrt{\Delta}/11$ neighbours with the same sum modulo $B$ as itself in $G$.
Note it is feasible, as otherwise (even taking into account the fact that $v_{n-1}$ and $v_n$ may always have the same sums modulo $B$),
$v_{n-1}$ and $v_n$ would have to have together at least $(2\sqrt{\Delta}/11-1)(B-12\sqrt{\Delta}-2\sqrt[3]{\Delta}+1)>2\Delta$ neighbours, a contradiction.
Then for the sum of $v_{n-1}$ we have at least $d_H(v_{n-1})-1$
available options (which form an arithmetic progression with step size $B$)
via
adjustments stemming from (i)
on the edges incident with $v_{n-1}$ in $H$ except on $v_{n-1}v_n$ and one more of these edges (which shall be specified below),
hence at least $3\sqrt{\Delta}/16-\sqrt[3]{\Delta}-1$ options due to (II).
These options thus contain at least
$(3\sqrt{\Delta}/16-\sqrt[3]{\Delta}-3)/2>\sqrt{\Delta}/11$
(disjoint) pairs from $\mathcal{P}$, at least one of which, say $P'$ is thus (due to our choice of the colour for $v_{n-1}v_n$) assigned to at most one neighbour of $v_{n-1}$ other than $v_n$ in $G$. Denote it by $v_k$ if it exists ($k<n-1$) -- then $v_{n-1}v_k$ shall be an edge incident with $v_{n-1}$ which we shall not modify.
We then perform the admissible modifications of the colours of the edges incident with $v_{n-1}$ so that $d^t_{c_T}(v_{n-1})\in P'$ and $d^t_{c_T}(v_k)\neq d^t_{c_T}(v_{n-1})$ (if $v_k$ exists). Set $P_{v_{n-1}}=P'$.
Then analogously for the sum of $v_{n}$ we have at lest $3\sqrt{\Delta}/16-\sqrt[3]{\Delta}-2$ available options (which form an arithmetic progression with step size $B$) via
admissible
adjustments
on the edges incident with $v_{n}$ except on $v_{n-1}v_n$, $v_{n}v_k$ (if it exists) and at most one additional one,
discussed below.
These options thus contain more than $\sqrt{\Delta}/11$
(disjoint) pairs from $\mathcal{P}$, at least one of which, say $P''$ is thus assigned to at most one neighbour of $v_{n}$ in $G$. Denote it by $v_{k'}$ if it exists ($k'<n$) -- then $v_{n}v_{k'}$ shall be an edge which we shall not modify.
We then perform the admissible modifications of the colours of the edges incident with $v_n$ so that $d^t_{c_T}(v_{n})\in P''$ and $d^t_{c_T}(v_{k'})\neq d^t_{c_T}(v_{n})$. Set $P_{v_{n}}=P''$.
Finally observe that in cases when $v_{n-1}\notin V_L$ or $v_{n}\notin V_L$ we may use the same strategy as above to obtain our goal, just skipping the unnecessary parts.
Note also that after analyzing $v_n$, the sums of the vertices in $H_1,\ldots,H_q$ shall not be modified in the further part of the algorithm.

At the end of the algorithm we thus have $d^t_{c_T}(u)\neq d^t_{c_T}(u)$ for every $uv\in E$ with $u,v\in V_L$, as claimed.

\subsection{Final Adjustment for Small Degree Vertices}

After analyzing all components of $H$, we greedily, one by one analyze all edges in $G_S$, modifying them if necessary so that there are no sum conflicts between
any vertex $v$ in $V_S$ with $d_{G_S}(v)\geq 1$ and its neighbours. We can do this, since we have more than $4(\Delta/4-1)+2$ colours available.
Note that afterwards all neighbours are sum-distinguished in $G$ as due to our construction if a vertex $u$ was in conflict with $w$ in $G$, then one of these two vertices, say $u$, would have to be an isolated vertex of $G_S$ which does not belong to $V_L$, while the other one -- $w$ -- would belong in $V\smallsetminus V_S$. Then however, $d_G(u)\leq \Delta/32$ and $d_G(w)\geq \Delta/4$, and as our colouring is proper and (\ref{maximal_colour_ineq}) holds, it is easy to check that
$$d^t_{c_T}(w)\geq \frac{\Delta^2}{32}>\frac{\Delta^2}{32}-\frac{\Delta^2}{2\cdot 32^2}+O(\Delta^\frac{3}{2})\geq d^t_{c_T}(u).$$
Note that we have not changed vertex colours within our construction, hence starting with a proper edge colouring yields a proper edge colouring on output.

\section{Remarks}
Note we did not strive to optimize the constant $95$ in our main result above.
This can be improved. It is in particular fairly easy in case of regular graphs.
Also when one considers the total and edge cases separately this constant can be reduced.

\end{document}